\documentclass[8pt]{article}
\usepackage[utf8]{inputenc}
\usepackage[english]{babel}
\usepackage[margin=1.2cm]{geometry}
\usepackage{mathtools}
\usepackage{wrapfig}
\usepackage{float}
\usepackage{pgf,tikz,pgfplots}
\pgfplotsset{compat=1.15}
\usepackage{mathrsfs}
\usetikzlibrary{arrows}
\usetikzlibrary{patterns}
\usepackage{xcolor}
\usepackage{iflang}

\usepackage{xcolor}

\usepackage{amsmath,amsthm,amssymb}
\usepackage{hyperref}
\usepackage{enumerate}

\usepackage{amsthm} 
\newtheorem{theorem}{Theorem}[section]

\newtheorem{question}[theorem]{Question}

\theoremstyle{remark}

\newcommand{\St}{{\mathcal St}}

        \newcommand{\HH}{{\mathcal H}}
        \renewcommand{\H}{\HH^1}

        \newcommand{\forget}[1]{}

\title{Open problems on Steiner trees and maximal distance minimizers}
\author{Yana Teplitskaya}
\date{Laboratoire de Math{\'e}matiques d’Orsay, Universit{\'e} Paris-Saclay, CNRS, Orsay, France}

\begin{document}

\setcounter{page}{1}

%  ----------------------------------------------

\title{Open problems on Steiner trees and maximal distance minimizers}
%

%\author{Danila Cherkashin}

%\affil{Institute of Mathematics and Informatics, Bulgarian Academy of Sciences, Sofia}

%  ----------------------------------------------

\maketitle

%  ----------------------------------------------

%\ReceivedAccepted{}{}

%  ----------------------------------------------

\abstract{In this work, I collect and discuss a series of open questions in one-dimensional geometric optimization in Euclidean spaces. The focus is on two classes of problems: maximal distance minimizers and Steiner trees. Maximal distance minimizers concern  finding a connected set of minimal length whose closed $r$-neighborhood covers a given compact set, whereas Steiner trees aim to find a minimal-length set connecting a prescribed set of points. For both problems, I briefly summarize known results and highlight the remaining open questions. While some questions can be approached with elementary methods, others remain highly challenging.}

%  ----------------------------------------------

\vskip-2mm
%\keywords {Steiner trees, maximal distance minimizers.}

%  ----------------------------------------------

%\AMSsc{2020}{49Q10, 49Q20, 05C05, 49Q22, 90B06}

%  ----------------------------------------------

\maketitle

%\section*{Introduction}

%\textcolor{red}{make more formal}
%While some questions can be approached with elementary methods, others remain highly challenging even for experts.

%In this work I tried to collect problems that are interesting to me, which I either actively think about or have already despaired of thinking about.
%First of all, these are problems of optimizing one-dimensional shapes.
%To understand this text, no knowledge is required; a junior student, as well as an experienced mathematician, can simply read the definitions and immediately begin solving non-trivial open problems.

%Please feel free to reach me out to discuss these problems, ask questions, share ideas, or collaborate on any of them.

\section{Maximal distance minimizers}
%\begin{problem}
For a given compact set $M \subset \mathbb{R}^d$ and a given positive number $r > 0$, consider the problem of finding a connected compact set $\Sigma$ of minimal length (one-dimensional Hausdorff measure $\mathcal{H}^1$) such that 
\[
M \subset \overline{B_r(\Sigma)},
\] 
where $\overline{B_r(\Sigma)}$ denotes the closed $r$-neighborhood of the set $\Sigma$.  
We refer to such a set $\Sigma$ as a  
\textit{Maximal distance minimizer} (\textit{MDM} or $r$-minimizer). %Existence of a solution is guaranteed (see~\cite{PaoSte04max}). 
A solution always exists and has no loops (homeomorphic image of a circle), see~\cite{paolini2004qualitative}.
For further details on MDMs, the reader is referred to our survey~\cite{cherkashin2022overview}.
%For a given compact set $M \subset \R^d$ and a given positive number $r > 0$ to find a connected compact set $\Sigma$ of the minimal length (one-dimensional Hausdorff measure $\H$) such that $M \subset \overline{B_r(\Sigma)}$, where $\overline{B_r(\Sigma)}$ denotes closed $r$-neighborhood of desired set $\Sigma$.  
%We will call such a solution $\Sigma$ a
% 

%A solution always exists and has no loops (homeomorphic image of a circle), see~\cite{paolini2004qualitative}.
%For detailed information on MDM we refer the reader to our review~\cite{cherkashin2022overview}.

 %\caption{A horseshoe}
 %   \label{horseshoe}

\subsection{Explicit examples}
\begin{question}
    What is a maximal distance minimizer for $M=\partial B_R$ with $R>r$? The problem remains open for $4.98r>R>r$.
\end{question}

\begin{figure}
  \begin{center}
  		\definecolor{yqqqyq}{rgb}{0.,0.,1.}%{0.5019607843137255,0.,0.5019607843137255}
		\definecolor{xdxdff}{rgb}{0.,0.,0.}
		\definecolor{ffqqqq}{rgb}{0.,0.,1.}
    \definecolor{ttzzqq}{rgb}{0.,0.,0.}
    \begin{tikzpicture}[line cap=round,line join=round,>=triangle 45,x=1.0cm,y=1.0cm]
    \clip(-5.754809259201689,-2.15198488708635626) rectangle (0.492320066152345,2.2268250892462824);
    \draw [very thick, color=ttzzqq] (-3.,2.)-- (-2.,2.);
    \draw [very thick, color=ttzzqq] (-2.,-2.)-- (-3.,-2.);
    \draw [very thick, color=ffqqqq] (-3.,-1.)-- (-2.,-1.);
    \draw [dotted] (-3.,1.)-- (-2.,1.);
    \draw [shift={(-3.,0.)},very thick, color=ttzzqq]  plot[domain=1.5707963267948966:4.71238898038469,variable=\t]({1.*2.*cos(\t r)+0.*2.*sin(\t r)},{0.*2.*cos(\t r)+1.*2.*sin(\t r)});
    \draw [shift={(-3.,0.)},dotted]  plot[domain=1.5707963267948966:4.71238898038469,variable=\t]({1.*1.*cos(\t r)+0.*1.*sin(\t r)},{0.*1.*cos(\t r)+1.*1.*sin(\t r)});
    \draw [shift={(-2.,0.)},dotted]  plot[domain=-1.5707963267948966:1.5707963267948966,variable=\t]({1.*1.*cos(\t r)+0.*1.*sin(\t r)},{0.*1.*cos(\t r)+1.*1.*sin(\t r)});
    \draw [shift={(-2.,0.)},very thick,color=ttzzqq]  plot[domain=-1.5707963267948966:1.5707963267948966,variable=\t]({1.*2.*cos(\t r)+0.*2.*sin(\t r)},{0.*2.*cos(\t r)+1.*2.*sin(\t r)});
    \draw [dash pattern=on 2pt off 2pt] (-2.5008100281931047,2.) circle (1.cm);
    \draw [very thick,color=yqqqyq] (-3.7309420990521054,0.6824394829091469)-- (-3.1832495111022516,1.2690579009478955);
    \draw [very thick,color=yqqqyq] (-1.8178907146013357,1.2695061868000874)-- (-1.2695061868000874,0.6829193135917684);
    \draw [shift={(-2.,0.)},very thick,color=ffqqqq]  plot[domain=-1.5707963267948966:0.7517515639553677,variable=\t]({1.*1.*cos(\t r)+0.*1.*sin(\t r)},{0.*1.*cos(\t r)+1.*1.*sin(\t r)});
    \draw [shift={(-3.,0.)},very thick,color=ffqqqq]  plot[domain=2.390497746093128:4.771687465439039,variable=\t]({1.*1.*cos(\t r)+0.*1.*sin(\t r)},{0.*1.*cos(\t r)+1.*1.*sin(\t r)});
    \begin{scriptsize}
    \draw[color=ttzzqq] (-4.8089822121844977,0.2025724503290166) node {$M$};
    %\draw [fill=xdxdff] (-2.5008100281931047,2.) circle (1.5pt);
    %\draw [fill=ffqqqq] (-3.1832495111022516,1.2690579009478955) circle (1.5pt);
    %\draw [fill=ffqqqq] (-1.8178907146013357,1.2695061868000874) circle (1.5pt);
    \draw[] (-3.763965847825129,-0.1818662773850787) node {$\Sigma$};
    \draw[color=black] (-2.059620164509451,0.8028001329511188) node {$M_r$};
    \end{scriptsize}
    \end{tikzpicture}
  \caption{A horseshoe for closed curve with big radius of curvature}
  \label{fig:horseshoe}
  \end{center}
\end{figure}
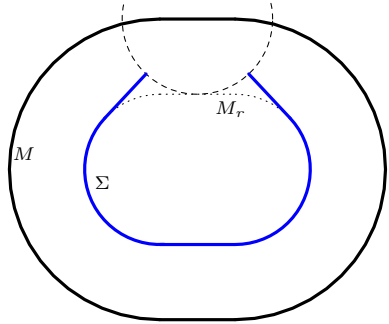

%\input{pic/horseshoe}
%\textbf{Partial answer and conjecture.}
Miranda, Paolini and Stepanov conjectured in~\cite{miranda2006one} that any minimizer for a circle of radius $R > r$ is a horseshoe id est a union of an arc of a circle of radius $R-r$ and two segments of its tangent lines.%\textcolor{red}{define horseshoe, $M_r$} 

The conjecture is confirmed in~\cite{cherkashin2018horseshoe} for $R>4.98r$. For smaller values of $4.98r>R > r$ the hypothesis remains unproven.

\begin{question}
     What is the set of maximal distance minimizers for $R$-stadium (which is the boundary of $R$-neighborhood of a segment)? Open for $r<R<5r$.
\end{question}

%\textbf{Partial answer and some information}
 %Let us call the boundary of an $R$-neighborhood of an interval a \textit{$R$-stadium}. It is obvious that the $R$-stadium has a radius of curvature of at least $R$ at each point. For $R < 1.75r$, for a sufficiently long stadium, the horseshoe is not an $r$-minimizer, since there is a set $\Sigma'$ that has a shorter length than any horseshoe and satisfies the other requirements.
In~\cite{cherkashin2018horseshoe} it is proved that for a closed curve with minimal radius of curvature $R>5r$ any MDM is a horseshoe, that is a union of an arc of a curve parallel to an arc of $M$ (and hence having a radius of curvature at least $R-r$) and two segments of its tangent lines. Thus, for $R>5r$, every solution for $R$-stadium is a horseshoe. It is also known (see~\cite{cherkashin2022overview}) that for $R<1.75r$ the solution is not a horseshoe. 

\begin{figure}
  \begin{center}
  %\begin{figure}[h]
 %   \centering

\begin{tikzpicture}[scale=1.3]

\foreach \y in {0, 3.51389977019889} {
    \draw[very thick, loosely dotted] 
        (-1, \y) -- (0, \y);
    \draw[very thick]
        (0, \y) -- (3.66, \y);
    \draw[very thick, loosely dotted] 
        (3.66, \y) -- (4.66, \y);
}
\draw[very thick] (-1, 0) arc (270:90:1.757);
\draw[very thick] (4.66, 0) arc (-90:90:1.757);

%\draw[thin, dashed] (0.2614137092297641 + 1, 0.0) arc (0:180:1);
%\draw[thin, dashed] (0.7842411276892922 + 1, 0.0) arc (0:180:1);
%\draw[thin, dashed] (1.3070685461488205 + 1, 0.0) arc (0:180:1);
\draw[thin, dashed] (1.8298959646083486 + 1, 0.0) arc (0:180:1);
%\draw[thin, dashed] (2.3527233830678767 + 1, 0.0) arc (0:180:1);
%\draw[thin, dashed] (2.875550801527405 + 1, 0.0) arc (0:180:1);
%\draw[thin, dashed] (3.398378219986933 + 1, 0.0) arc (0:180:1);
%\draw[thin, dashed] (1.8298959646083486 + 1, 3.75) arc (0:-180:1);

\draw[blue, thick] 
(0.0, 0.9652268503449002) -- 
(1.5684822553785844, 0.9652268503449002) --
(1.8298959646083486, 1.11615412573856) --
(2.0913096738381127, 0.9652268503449002) -- 
(3.659791929216697, 0.9652268503449002);

\node[blue] at (1.61,1.8) {$\Sigma_0$};

\draw[dashed] (0.0, 0.0) -- (0.0, 3.51389977019889);
\draw[dashed] (3.659791929216697, 0.0) -- (3.659791929216697, 3.51389977019889);

\fill[blue] (0.0, 0.9652268503449002) 
    circle (1.5pt) node[black, below right] {$a$};
\fill[blue] (3.659791929216697, 0.9652268503449002) 
    circle (1.5pt) node[black, below left] {$b$};

\draw[blue, thick] 
(1.8298959646083486, 1.11615412573856) -- 
(1.8298959646083486, 2.58208251632886);

\draw[thin, dashed] (1, 3.51389977019889) arc (0:-90:1);
\draw[thin, dashed] (2.66, 3.51389977019889) arc (180:270:1);

\draw[blue, thick]
(1.8298959646083486, 2.58208251632886) -- 
(0.9149479823041743, 3.11032798020668);

\draw[blue, thick] 
(1.8298959646083486, 2.58208251632886) -- 
(2.744843946912523, 3.11032798020668);

\end{tikzpicture}

   % \caption{Horseshoe is not a minimizer for long enough stadium with $R < 1.75r$.}
  %  \label{stadion}
%\end{figure}
  \caption{A competitor shorter than the horseshoe for a stadium} 
  \label{fig:stadium}
  \end{center}
\end{figure}
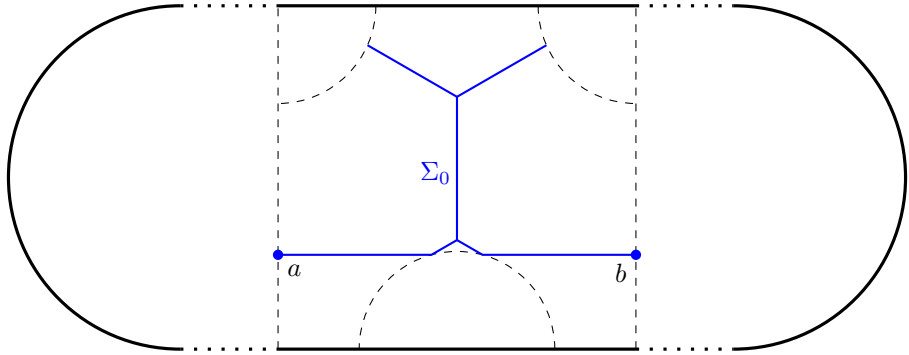

%\begin{figure}{r}{0.1\textwidth}
 %  \vspace{-1cm}
  %\begin{center}
  %\input{pic/stadium2}
  %\end{center}
  %\caption{A horseshoe}
  %\label{fig:horseshoe}
%\end{figure}
 Thus the following sub-questions remain open:
\begin{itemize}
    \item To find the supremum of $C$ such that there exists MDM for $Cr$-neighborhood which is not horseshoe.
    \item To find the infimum of $c$ such that for any  closed curve with minimal radius of curvature $cr$ every MDM is horseshoe.
\end{itemize}
It might also be helpful to know that any minimizer for $R$-stadium consists of line segments and curves parallel to the curves contained in $R$-stadium (these curves coincide with arcs of circles of radius $R-r$).

\begin{question}
  What is the set of MDMs for any explicit set you like?
 \end{question}

  For example for 
    \begin{itemize}
        \item A disc or ball $B_R$ in $\mathbb{R}^d$, $d\geq 2$. 
        Very challenging question, nothing is known except the usual regularity properties such as absence of loops.
        \item Rectangle (the boundary). For a small enough $r>0$ the question is answered, see~\cite{cherkashin2021maximal}.
        \item (Regular) $n$-gons.  
        Looks like for $n\geq 6$ MDM for regular $n$-gon will be just a poly-line.  Although this question seems less difficult than finding minimizers for other sets, no one has written a clean proof yet. By the regularity properties it is clear that the solution is a finite union of line segments with angles at least $2\pi/3$.
        The question is interesting already for $n=3$. Our methods for rectangle might help, see~\cite{cherkashin2021maximal}. %\textcolor{red}{$r$ should be small otherwise something awful happens}
    \end{itemize} 

\subsection{Regularity and uniqueness}
\begin{question}
Is any MDM in $\mathbb{R}^d$ a union of a finite number of injective curves? Open for $d\geq 3$. 
\end{question}

%\textbf{Some information}.
In $\mathbb{R}^2$ the answer is positive, see~\cite{gordeev2022regularity}.
In higher dimensions some information about local behavior is proved in the same paper: there are at most three tangent rays at any point of minimizers and the angles between the tangent rays are at least $2\pi/3$.
Recall that a solution always exists and never contains loops, see~\cite{paolini2004qualitative}.
%, can be not unique but usually is (see~\cite{basok2018uniqueness}).

\begin{question}
Let $\Sigma$ be an $r$-minimizer for some $M$. Is $\Sigma$ the unique $r$-minimizer for $\overline{B_r(\Sigma)}$?
\end{question}
Clearly that if $\Sigma$ is an $r$-minimizer for some $M$, then it is a minimizer for $\overline{B_r(\Sigma)}$, see~\cite{basok2022inverse}. A weaker form of this question is if we replace $r$ with some positive $r_0 < r$ in the hypothesis.

Some points of an MDM are more important than others: we call a point $x\in \Sigma$ \textit{energetic} if there exists a point $y \in M$ such that $|xy|=r$ and $B_r(y) \cap \Sigma = \emptyset$. We denote by $G_\Sigma$ the set of energetic points of $\Sigma$.
%for every $\varepsilon>0$

\begin{question}
    Fix $M \subset \mathbb{R}^d$ and $r > 0$.
Does the set $G_\Sigma$ of energetic points determine a minimizer $\Sigma$?
\end{question}

%\textcolor{red}{maybe to add that the connection between questions 2.5 and 2.6 is a nice question itself}\\

The connection between Questions~1.5 and 1.6 is a nice question itself.

\subsection{Algebraic point of view}

Reflections from a spherical mirror play an important role in geometry of MDMs. In particular MDM for a $3$-point set $M$ often coincides with a solution to the Alhazen’s billiard problem, 
which is actually reduced to solving a certain fourth-degree equation and thus cannot be solved by a straightedge (ruler) and a compass, see~\cite{fujimura2019ptolemy}. 
The following question is open even for $n = 4$. 

\begin{question}
Let $M \subset \mathbb{R}^2$ be a finite set of cardinality $n$. 
How does the explicit description of a minimizer $\Sigma$ depend on $M$ in terms of the degrees of the resulting equations?
\end{question}

For simplicity one may assume that the points from $M$ have rational coordinates. 
Then it is enough to find the largest degree of the extension of $\mathbb{Q}$ generated by the set of energetic points $G_\Sigma$.

%Problem can be stated in any dimension.% but even at the plane it seems not to be elementary. 
%\begin{itemize}
      
 %    \item In this case, the statement of the Theorem~\ref{horseshoeT} without the condition on the minimum radius of curvature becomes incorrect, as we will show below.

%For $i \in \N$ we define the set $\Sigma_i$ as a local Steiner tree, shown in Fig.~\ref{stadion}.
%Let $\Sigma'$ consist of several sets $\Sigma_i$ located along the entire length of the stadium, so that $\Sigma_i$ and $\Sigma_{i+1}$ are glued together at points $b_i$ and $a_{i+1}$.
%For $R < 1.75r$, the length of $\Sigma_i$ is strictly less than $2|ab|$. Thus, for a sufficiently long stadium, $\Sigma'$ has length $cL + O(1)$, where $L$ is the length of the stadium, and $c < 2$ is a constant depending on $\Sigma_1$ and $R$. At the same time, the horseshoe has length $2L + O(1)$.

%This example motivates the next two questions.

 % \begin{question}
  %For what minimum $c$ does the Theorem~\ref{horseshoeT} hold with $5r$ replaced by $cr$ in the condition?
  %\end{question}
 
  %\begin{question}
%What do a lot of stadium minimizers look like?
 % \end{question}

 \section{Steiner trees}

    For a given compact set $\mathcal{A} \subset \mathbb{R}^d$ find a closed set $\St$ with minimal length (one-dimensional Hausdorff measure $\H$) such that $\St \cup \mathcal{A}$ is connected.
%If $\mathcal{A}$ is totally disconnected, $\St$ is connected
We will call such a solution $\St$ \textit{Steiner set}. If $\St$ is connected (which is the case if $\mathcal{A}$ is finite or totally disconnected) we call it Steiner tree. In 2013, Paolini and Stepanov showed the existence of a solution, see~\cite{paolini2013existence}. If a solution $\St$ has finite length, then it is compact, has no loops, and consists of countably many connected components.
%The points of $\mathcal A$ belonging to the closure of $\St$ will be called \textit{terminals}. 
The edges of the locally finite graph $\St \setminus \mathcal A$ are straight line segments.
The maximal degree (in graph-theoretic sense) of a vertex is at most $3$. Moreover, only terminals can have degree $1$ or $2$, all the other vertices have degree $3$ and are called \textit{Steiner} (branching) points. The angle between any two adjacent edges of a $\St$ is at least $2\pi/3$, in particular $\St$ in a neighborhood of a branching point is a regular tripod (all three angles are equal to $2\pi/3$) and hence are coplanar. If every terminal has degree 1 then $\St$ is called \textit{full}.

%\textit{different}
\subsection{Uniqueness and existence}
\begin{question}
    Does there exist a finite set of terminals in $\mathbb{R}^d$, $d \geq 3$ with two different Steiner trees $\St_1$ and $\St_2$ such that the directions of all segments incident to terminals are the same (in other words trees coincide in a neighborhood of terminals)? 
\end{question}
In the plane K. Oblakov has proved that it is impossible~\cite{oblakov2009non}.% Proof is nice but very planar. 

%, so I don't know how to generalize it for higher demensions.
%\textcolor{red}{maybe place to say about edelsbruner and topology}
% every finite planar input the directions of all segments incident to terminals determine ST\cite{oblakov2009non}

\begin{question}   \label{questUniq}
Is a solution to the Steiner tree problem unique for a general input? 
\end{question}

We~\cite{basok2018uniqueness} provided the positive answer for finite planar input. 
Moreover we proved that in $\mathbb{R}^2$ the set of $n$-point ambiguous configurations has the
Hausdorff dimension at most $2n-1$ (where a configuration is called \textit{ambiguous} if there are several Steiner trees for it). This particularly means that $n$-dimensional Hausdorff measure of ambiguous $n$-point configurations is $0$, id est a general planar $n$-input has a unique Steiner tree.

\begin{question}%[\cite{paolini2012minimal}]
Is it possible to find a bounded set $A \subset \mathbb{R}^d$ such that the problem:
``find $\mathcal{S}$ such that $\H(\mathcal{S})$ is minimal among all sets for which $\mathcal{S} \cup A$ is connected''
does not have a solution?
\end{question}
This question was asked by E. Paolini, who is interested whether the assumption that $A \subset \mathbb{R}^d$ is closed is necessary for the problem to be solvable.   

\subsection{Gilbert--Pollak conjecture on the Steiner ratio}
\textit{A Minimum spanning tree (MST)} is a minimal connection of a given set $A$ of points in the plane by segments with endpoints in $A$.
 \textit{The Steiner ratio} is the infimum of the ratio of the total length of Steiner tree to the minimum spanning tree  for a finite set of points in the Euclidean plane.
\begin{question} 
What is the Steiner ratio in $\mathbb{R}^d$, $d\geq 2$?

%Where \textit{the Steiner ratio} is the infimum of the ratio of the total length of Steiner tree to the minimum spanning tree  for a finite set of points and 
%\textit{a Minimum spanning tree} for $A$ is a minimal connection of a given set $A$ of points by segments with endpoints in $A$.
\end{question}

%It is well known that an MST can be found in a polynomial time in the size of input (via Prim's or Kruskal's algorithm).

%\textbf{Conjectures}
This is the most well-known open problem on Steiner tree, so called \textit{Gilbert--Pollak conjecture on Steiner ratio}, id est whether the planar Steiner ratio is reached on an equilateral triangle (and so is equal to $\frac{\sqrt{3}}{2}=0.866\ldots$). This conjecture is still open (as of 2025 year), see~\cite{ivanov2012steiner}. Graham~\cite{graham2019some} offered 1000\$ for the proof (in his work, the reader can find the main information about the problem and attempts to solve it). The best proven bound is about $0.842$ (1985,~\cite{chung1985new}).

In higher dimensions the problem becomes even more difficult (but cheaper): 500\$ prize is offered for a conjecture in $\mathbb{R}^3$, higher dimensions are priceless. 

As Ronald L. Graham and Frank K. Hwang showed in 1976 in any Euclidean space, the Steiner ratio is at most $\sqrt{3}$, see~\cite{graham1976remarks}.
It is also known that for any metric space, the Steiner ratio is at most two.

The ratio for the vertices of a regular simplex, including the tetrahedron, is better than the planar value, but it is not optimal. D. Smith and J. MacGregor Smith constructed in~\cite{smith1995steiner} so called $d$-sausage configurations which provide the best known and conjecturally optimal value of the Steiner ratio in $\mathbb{R}^d$, $d\geq 2$. For $d=3$ this value is approximately $0.784$.

\subsection{Density type questions}
These questions make sense in both two-dimensional Euclidean space and in a higher dimension. They were suggested to me in personal conversation with Guy David. % I don't know what is the best statement of the problem.
%\textcolor{red}{rewrite everything here :( }

Let $\St$ be a Steiner tree and suppose that $B_r(x)$ contains no terminal points i.e. $\mathcal{A} \cap B_r(x))=\emptyset$. We will additionally assume that the set has a finite length:
\[
\H(\Sigma \cap \overline{B_{r}(x)}) <\infty.
\]
In the plane this condition could be replaced by $\mathcal{A}\subset \partial B_r(x)$. We can also add requirements such as that $\St \cap \overline{B_r}$ is connected or full.

%And ask the questions like
We are interested in questions similar to the following.
\begin{question}
    What is the largest possible number of $\sharp(\partial B_{tr}(x) \cap \Sigma )$ for $t=\frac{1}{2}$ (or other $t \in (0;1)$)? 
\end{question}

    \begin{figure}
  \begin{center}
  %\begin{figure}[h]
 %   \begin{center}
     \definecolor{qqwuqq}{rgb}{0.,0.39215686274509803,0.}
\definecolor{ffvvqq}{rgb}{1.,0.3333333333333333,0.}
\definecolor{uuuuuu}{rgb}{0.26666666666666666,0.26666666666666666,0.26666666666666666}
\definecolor{ududff}{rgb}{0.30196078431372547,0.30196078431372547,1.}
\begin{tikzpicture}[line cap=round,line join=round,>=triangle 45,x=1.0cm,y=1.0cm]
\clip(-2.610277148878922,-2.08487906707841) rectangle (3.760094427223567,3.2482136534371877);
\draw [line width=1.2pt] (0.,0.) circle (2.cm);
\draw [line width=1.2pt] (0.,0.) circle (1.cm);
\draw  [dotted] [line width=0.8pt,color=qqwuqq] (-1.937607863137521,0.4956569062442783)-- (-0.02078935840890737,-0.018026693260688544);
\draw [dotted]  [line width=0.8pt,color=qqwuqq] (-0.02078935840890737,-0.018026693260688544)-- (0.49362601332867106,-1.9381262494907916);
\draw [dotted]  [line width=0.8pt,color=qqwuqq] (-0.02078935840890737,-0.018026693260688544)-- (0.7982789688621812,0.8009829274547386);
\draw [dotted]  [line width=0.8pt,color=qqwuqq] (0.7982789688621812,0.8009829274547386)-- (0.49362601332867106,1.9381262494907918);
\draw [dotted]  [line width=0.8pt,color=qqwuqq] (0.7982789688621812,0.8009829274547386)-- (1.937607863137521,0.4956569062442783);
\draw [ultra thick, blue] (-1.937607863137521,0.4956569062442783)-- (0.26533726635090726,1.0860188945264015);
\draw [ultra thick, blue] (0.26533726635090726,1.0860188945264015)-- (0.49362601332867106,1.9381262494907918);
\draw [ultra thick, blue] (0.26533726635090726,1.0860188945264015)-- (1.0844055936219945,0.26700927381097483);
\draw [ultra thick, blue] (1.0844055936219945,0.26700927381097483)-- (1.937607863137521,0.4956569062442783);
\draw [ultra thick, blue] (1.0844055936219945,0.26700927381097483)-- (0.49362601332867106,-1.9381262494907916);
\begin{scriptsize}
%\draw [fill=ududff] (0.49362601332867106,0.4956569062442783) circle (2.5pt);
\draw [fill=uuuuuu] (-1.937607863137521,0.4956569062442783) circle (2.0pt);
\draw [fill=uuuuuu] (0.49362601332867106,1.9381262494907918) circle (2.0pt);
\draw [fill=uuuuuu] (1.937607863137521,0.4956569062442783) circle (2.0pt);
\draw [fill=uuuuuu] (0.49362601332867106,-1.9381262494907916) circle (2.0pt);
\draw [fill=blue] (0.26533726635090726,1.0860188945264015) circle (2.0pt);
%\draw [fill=qqwuqq] (0.7982789688621812,0.8009829274547386) circle (2.0pt);
\draw [fill=blue] (1.0844055936219945,0.26700927381097483) circle (2.0pt);
%\draw [fill=qqwuqq] (-0.02078935840890737,-0.018026693260688544) circle (2.0pt);
\end{scriptsize}
\end{tikzpicture}
 %\end{center}
  %  \label{fig:6 points}
% \caption{Example of $6$ points at the half circle}
%\end{figure}
  \end{center}
  \caption{Example of $6$ points at the half-radius circle}
  \label{fig:6points}
\end{figure}
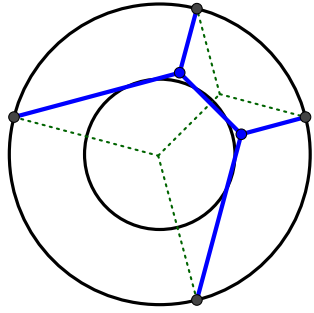
%\textbf{An example and conjecture} 
%For $t=1$ we need a requirement for $\Sigma$ to be full. 
For $t=\frac{1}{2}$ M. Khristoforov in personal communication provided a planar example for $6$ (see Fig.~\ref{fig:6points}). It is an open question if $6$ is the largest possible number at the plane.
   % \input{pic/6points}

%\begin{figure}
    %\centering
    %\includegraphics[width=0.5\linewidth]{}
   % \caption{Example of $6$ points at the half circle}
%    \label{fig:6 points}
%\end{figure}
\begin{question}
What is the largest possible number of branching points of $\St$ in $B_{tr}(x)$ for $t=\frac{1}{2}$ or other $t \in (0;1)$?
\end{question} 
%\textcolor{red}{ картинка от христофорова}
%\begin{question}
 %   Additionally assume that $\Sigma$ is full. What is the maximum of $\H(\St \cap B_{r}(x))/r$? 
%\end{question}
%\textbf{No limits close to the boundary}
 In work in progress Cherkashin and Prozorov show that the desired value is at most $\left( \frac{32 d}{1-t} \right)^{d-1}$. Also they provide a planar example of full connected Steiner tree with four limit points at the boundary $\partial B_r(x)$. Thus for $t$ goes to $1$ the values in both Questions 2.5 and 2.6 go to infinity.

\begin{question}
    Let $\St \subset \mathbb{R}^d$ for $d > 2$. Additionally assume that $\St$ is full. What is the maximum of $\frac{\H(\St \cap B_{tr}(x))}{r}$ for $t \in (0,1)$?
    %  for $t=\frac{1}{2}$
\end{question}
%\textbf{Upper bound estimation}
In work in progress, we (Cherkashin--Prozorov--T.) show that 
    \[
    \frac{\H(\St \cap B_{t r}(x))}{r} \leq \frac{(32d)^{d-2}}{(1-t)^{d-2}}.
    \]

 %$\St_\varepsilon$ is Ahlfors--David regular with the  parameters which depend only on $d$ and $\varepsilon$ (and do not depend on $\mathcal{A}$). Namely, for every $\varepsilon > 0$, $\rho \in (0,1)$ and $x \in \St_\varepsilon$ we have
  %  \[
   % \frac{\H(\St_\varepsilon \cap B_{\rho \varepsilon}(x))}{\varepsilon} \leq \left ( \frac{32d}{1-\rho} \right) ^{d-2}.
    %\]

%\begin{question}
%What is the maximum of $\H(\St \cap B_r(x))/r$?  
%\end{question}

%Let $\mathcal{A}$ be a finite subset of $\partial B_1(0)$ (or another good set). What is the largest possible $\sharp(\St \cap \partial B_{1/2}(0))$? 
%What is the largest possible number of branching points of $\St$ lying within $$B_{1/2}(0)$?
%What is the upper bound for $\H(\St \cap B_{1/2}(0))$?
%Similar questions can be asked for any radius in place of $\frac 1 2$.

\subsection{Connection with maximal distance minimizers}

A natural way to investigate the connection is to consider the case of a finite set $M$ in maximal distance minimizers problem. A minimizer for a finite set $M \subset \mathbb{R}^d$ is a finite union of at most $2\sharp M - 3$ segments. In this case the maximal distance minimizers problem comes down to connecting $r$-neighborhoods of all the points from $M$. 

A \textit{topology} $T$ of a labeled Steiner tree $\St$ is the corresponding abstract graph with labeled terminals and unlabeled Steiner points.
We showed that if a $\St$ for $n$ terminals in $\mathbb{R}^d$ preserves the topology for an arbitrary small enough perturbation of the terminals, then $\St$ is an $r$-minimizer for some $n$-tuple $M$ and a small enough $r$, see~\cite{basok2022inverse}.
Recall that in the plane a Steiner tree for a random input is unique with unit probability. Also in the plane a somehow opposite statement holds:
if $\St$ is a full planar Steiner tree for $n$ terminals which is not unique, then $\St$ cannot be a minimizer for any $M$ being an $n$-tuple of points. 

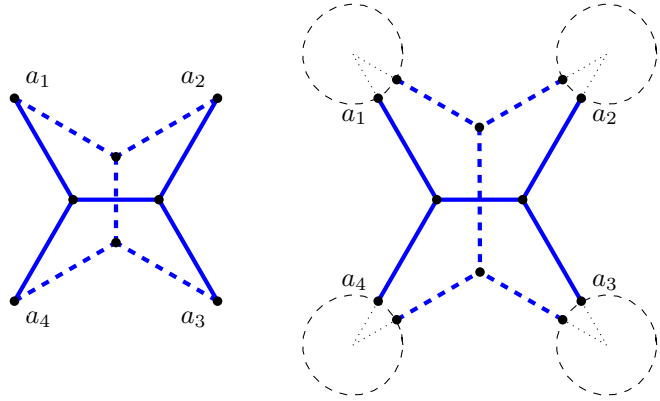
\begin{figure}
  \begin{center}
     \begin{tikzpicture}[scale=0.9]
    \def\r{1.5cm}
    \draw[ultra thick, blue]
        (-\r, \r) coordinate(x1) node[black, above right]{$a_1$} --++ (-60:{\r/cos(30)}) coordinate (x5);
    \draw[ultra thick, blue]
        (\r,\r) coordinate(x2) node[black, above left]{$a_2$} --++ (-120:{\r/cos(30)}) coordinate (x6);
    \draw[ultra thick, blue]
        (\r, -\r) coordinate(x3) node[black, below left]{$a_3$} --++ (120:{\r/cos(30)});
    \draw[ultra thick, blue]
        (-\r,-\r) coordinate(x4) node[black, below right]{$a_4$} --++ (60:{\r/cos(30)});
    \draw[ultra thick, blue]
        (x5) -- (x6);
        
    \draw [dotted,white] (x1)  --++ (-60:{-\r/2}) coordinate (x11);
    \draw [dotted,white] (x2)  --++ (-120:{-\r/2}) coordinate (x12);
    \draw [dotted,white] (x3)  --++ (120:{-\r/2}) coordinate (x13);
    \draw [dotted,white] (x4)  --++ (60:{-\r/2}) coordinate (x14);
    
    \draw [dashed,shift={(x11)},white]  plot[domain=0:6.5,variable=\t]({cos(\t r)/1.35},{sin(\t r)/1.35});
    \draw [dashed,shift={(x12)},white]  plot[domain=0:6.5,variable=\t]({cos(\t r)/1.35},{sin(\t r)/1.35});
    \draw [dashed,shift={(x13)},white]  plot[domain=0:6.5,variable=\t]({cos(\t r)/1.35},{sin(\t r)/1.35});
    \draw [dashed,shift={(x14)},white]  plot[domain=0:6.5,variable=\t]({cos(\t r)/1.35},{sin(\t r)/1.35});

    \draw[ultra thick, blue, dashed]
        (x1)  --++ (-30:{\r/cos(30)}) coordinate (x7);
    \draw[ultra thick, blue, dashed]
        (x2)  --++ (-150:{\r/cos(30)});
    \draw[ultra thick, blue, dashed]
        (x3)  --++ (150:{\r/cos(30)}) coordinate (x8);
    \draw[ultra thick, blue, dashed]
        (x4) --++ (30:{\r/cos(30)});
    \draw[ultra thick, blue, dashed]
        (x7) -- (x8);
    \foreach \x in{1,2,...,8}{
        \fill (x\x) circle (2pt);
    }
\end{tikzpicture}
\begin{tikzpicture}[scale=0.9]
    \def\r{1.5cm}
    \draw[ultra thick, blue]
        (-\r, \r) coordinate(x1) node[black, below left]{$a_1$} --++ (-60:{\r/cos(30)}) coordinate (x5);
    \draw[ultra thick, blue]
        (\r,\r) coordinate(x2) node[black, below right]{$a_2$} --++ (-120:{\r/cos(30)}) coordinate (x6);
    \draw[ultra thick, blue]
        (\r, -\r) coordinate(x3) node[black, above right]{$a_3$} --++ (120:{\r/cos(30)});
    \draw[ultra thick, blue]
        (-\r,-\r) coordinate(x4) node[black, above left]{$a_4$} --++ (60:{\r/cos(30)});
    \draw[ultra thick, blue]
        (x5) -- (x6);
        
    \draw [dotted] (x1)  --++ (-60:{-\r/2}) coordinate (x11);
    \draw [dotted] (x2)  --++ (-120:{-\r/2}) coordinate (x12);
    \draw [dotted] (x3)  --++ (120:{-\r/2}) coordinate (x13);
    \draw [dotted] (x4)  --++ (60:{-\r/2}) coordinate (x14);
    
    \draw [dashed,shift={(x11)}]  plot[domain=0:6.5,variable=\t]({cos(\t r)/1.35},{sin(\t r)/1.35});
    \draw [dashed,shift={(x12)}]  plot[domain=0:6.5,variable=\t]({cos(\t r)/1.35},{sin(\t r)/1.35});
    \draw [dashed,shift={(x13)}]  plot[domain=0:6.5,variable=\t]({cos(\t r)/1.35},{sin(\t r)/1.35});
    \draw [dashed,shift={(x14)}]  plot[domain=0:6.5,variable=\t]({cos(\t r)/1.35},{sin(\t r)/1.35});

    \draw [dotted] (x11)  --++ (-30:{\r/2}) coordinate (x21);
    \draw [dotted] (x12)  --++ (-150:{\r/2}) coordinate (x22);
    \draw [dotted] (x13)  --++ (150:{\r/2}) coordinate (x23);
    \draw [dotted] (x14)  --++ (30:{\r/2}) coordinate (x24);

    \draw[ultra thick, blue, dashed]
        (x21)  --++ (-30:{\r*0.94}) coordinate (x7);
    \draw[ultra thick, blue, dashed]
        (x22)  --++ (-150:{\r*0.94});
    \draw[ultra thick, blue, dashed]
        (x23)  --++ (150:{\r*0.94}) coordinate (x8);
    \draw[ultra thick, blue, dashed]
        (x24) --++ (30:{\r*0.94});
    \draw[ultra thick, blue, dashed]
        (x7) -- (x8);
    
    \foreach \x in{1,2,...,8}{
        \fill (x\x) circle (2pt);
        }
    \foreach \x in{1,2,...,4}{
        \fill (x2\x) circle (2pt);
    }
\end{tikzpicture}
  \end{center}
  \caption{A non-unique Steiner tree which is not a Maximal distance minimizer}
  \label{picnotunique}
\end{figure}

To illustrate this statement consider a square $a_1a_2a_3a_4$. There are two Steiner trees for $a_1,a_2,a_3,a_4$ (see the left-hand side of Fig.~\ref{picnotunique}), let us pick the solid one. The right-hand side of Fig.~\ref{picnotunique} shows that an $r$-minimizer for every positive $r$ has the topology of the dotted Steiner tree.

Thereby the picture is clear for a full Steiner tree for planar finite input. All possible generalizations remain open.
%(and quite important and interesting) questions:
\begin{question}
Let $\mathcal{A} \subset \mathbb{R}^d$ be a compact totally disconnected set with at least two Steiner trees, say $S$ and $T$. 
Can $S$ be a solution to the maximal distance minimizers problem?
\end{question}

The partial results are known in the case when $d=2$: if $S$ or $T$ is full then the answer is negative: 
a finite $\mathcal{A}$ is covered in~\cite{basok2022inverse}, and an infinite $\mathcal{A}$ contradicts decomposition of a maximal distance minimizers into a finite number of simple curves, see~\cite{gordeev2022regularity}.

\subsection{Large number of terminals and asymptotical behavior of Steiner trees}
Random geometric structures have been extensively studied~\cite{solomon1978geometric}~\cite{penrose2003random}. However, to the best of my knowledge, the Steiner tree problem (STP) has not yet been explored in a random setting. The closest related topics that have been investigated include minimum spanning trees~\cite{gilbert1965random}, the rectilinear STP~\cite{bern1988two} and the traveling salesperson problem~\cite{beardwood1959shortest}~\cite{frieze2016separating}. It is known that the asymptotics of the expected length of a solution for $N$ random points in the unit cube in $\mathbb{R}^d$ is $\beta(d) N^{\frac{d-1}{d}}$ and the constants $\beta(d)$ are distinct for different problems.

\begin{question} \label{questRandom}
Does the length of a Steiner tree for a random $N$-point input in $\Omega \subset \mathbb{R}^d$ have asymptotics $\beta N^{\frac{d-1}{d}}$, where $\beta$ is a constant depending on $\Omega$ and the random model?

%Let $\Omega \subset \mathbb{R}^d$ be a domain. To show that the length of a Steiner tree for a random $N$-point input has asymptotics $\beta N^{\frac{d-1}{d}}$,
%where $\beta$ is a constant, depending on $\Omega$ and the random model.
\end{question}

It seems highly ambitious to expect a complete answer to Question~\ref{questRandom} with an explicit $\beta$. Firstly, corresponding results are not yet known even for the better-studied problems mentioned earlier. Secondly, there is still a lack of comprehensive understanding of the Steiner tree problem for $d > 2$. Therefore, one may focus on establishing strong upper and lower bounds on the length of a random Steiner tree. 

\begin{question}
Can one obtain the bounds of form $\beta_1 N^{\frac{d-1}{d}} \leq \mathbb{E} \mathcal{H}^1(\St) \leq \beta_2 N^{\frac{d-1}{d}}$ that distinguish the Steiner tree problem from the other mentioned problems (the minimum spanning tree, the rectilinear Steiner tree and the traveling salesperson problems)? What are good upper and lower bounds for the Steiner ratio in a random setup?
        %minimum spanning trees~\cite{gilbert1965random}, the rectilinear STP~\cite{bern1988two} and the traveling salesperson proble
       % To obtain the bounds of type $\beta_1 N^{\frac{d-1}{d}} \leq \mathbb{E} \mathcal{H}^1(\St) \leq \beta_2 N^{\frac{d-1}{d}}$ that separate STP from the other mentioned problems (MST, RSP, TSP). Find a good upper and lower bound of the Steiner ratio in a random setup. Show that a solution to STP is unique for a countable input.
\end{question}

The next question is the result of personal conversations between G. Buttazzo, D. Cherkashin and E. Paolini.

Consider a planar connected compact set $\Omega$. 
Let $L_\Omega (n)$ be the maximal possible length of a Steiner tree for $n$-point set $\mathcal A \subset \Omega$. 
It is well-known that $L_\Omega (n)$ grows like $\sqrt{n}$. Although this is not written explicitly in the literature, standard methods should imply
%It's not written anywhere, but standard methods should give
the existence of limit
\[
l = \lim \sqrt{\frac{L_\Omega(n)}{\mathcal{H}^2(\Omega)}}.
\]

The (asymptotically) worst known set~\cite{chung1981steiner, few1955shortest} is obtained by the  intersection of the hexagonal lattice with $\Omega$, which gives
\[
l \geq \left(\frac{4}{3} \right)^{1/4} = 1.07\dots.
\]
This bound seems to be sharp.

To prove the reverse inequality E. Paolini suggested to prove that one can always connect $n$ points in a rectangle 
$\sqrt 3 \times (n-1)$ by a set of length at most $\sqrt{3}(n-1)+o(1)$. 
This bound is assumed to be achieved on the set of points with coordinates $(2i, 0)$ and $(2i+1, \sqrt{3})$.

Therefore we arrive at the following question. 

\begin{question}
     Is it true that one can always connect $n$ points in a rectangle 
$\sqrt 3 \times (n-1)$ by a set of length at most $\sqrt{3}(n-1)+o(1)$?
\end{question} 

%\textcolor{red}{to explain what is topology? to skip this question?}
%\begin{remark*}
    The same problems are also open for minimum spanning trees instead of Steiner trees. The extremal sets are conjectured to be the same. The asymptotic ratio of these extremal sets for Steiner and minimum spanning trees is $\frac{\sqrt{3}}{2}$ which coincides with the one in Gilbert--Pollak conjecture.
%\end{remark*} 

The last question of this paragraph was stated in our work~\cite{cherkashin2025steiner}. 
\begin{question}
    It is known, see~\cite{gilbert1968steiner}, that the total number of full Steiner topologies (the underlying abstract graphs) with $n$ terminals is 
\[
\frac{(2n-4)!}{2^{n-2}(n-2)!}.
\]
How many of them may simultaneously give a shortest solution to the Steiner problem? 
\end{question}
In all known examples the number is at most exponential in $n$ while the number of topologies has a factorial order of growth. Note that the answer may depend heavily on the dimension.

\subsection{Nonplanar explicit examples}
Recall that Steiner trees in a Euclidean $d$-dimensional space with $d \geq 3$ have the same local structure as in the plane, id est every branching point has degree $3$ and moreover three adjacent segments have pairwise angles $2\pi/3$ and thus are coplanar. However, a global structure is much more complicated. In particular, there is no explicit algorithm which may be caused by the fact that Steiner points are no longer expressed in radicals~\cite{smith1992find}.

% In particular a full tree may have no wind-rose.

To the best of our knowledge all known series of explicit solutions are planar. Even for the vertices of regular simplices~\cite{fleischmann2025steiner} in high dimensions the Steiner trees remain unknown.
% \begin{figure}{r}{0.3\textwidth}%{0.1\textwidth}
%   \vspace{-1cm}
%  \begin{left}
  %\input{pic/angletree.tex}
 % \end{left}
  %\caption{Steiner tree from~\cite{cherkashin2025steiner}}
 % \label{fig:angle}
%\end{figure}

\begin{figure}
  \begin{center}
  \begin{tikzpicture}
%\begin{scope}[shift={(15,0)}]
  
  \def\lambda{0.3}
  \def\ty{0.4} %% 1-2*\lambda
  \def\r{7cm}
  \def\rr{\lambda*\r}
  \def\rrr{\ty*\ty*\r}
  \def\rrrr{\ty*\ty*\ty*\r}
  \def\rrrrr{\ty*\ty*\ty*\ty*\r}

    \path (-1.5*\r,0) coordinate [label=above right:$y_0$] (y0);
    \path (-\r,0) coordinate [label=above left:$y_1$] (y1);

	\path (y1) ++(60:0.75*\ty*\r) coordinate [label=left:$y_2$] (y2);
	\path (y1) ++(-60:0.75*\ty*\r) coordinate [label=left:$y_3$] (y3);
	\path (y2) ++(120:\ty*\rr) coordinate [label=right:$y_4$] (y4);
	\path (y2) ++(0:\ty*\rr) coordinate [label={[xshift=-0.1cm, yshift=0.1cm]:$y_5$}] (y5);
	\path (y3) ++(0:\ty*\rr) coordinate [label={[xshift=-0.2cm, yshift=-0.42cm]:$y_6$}] (y6);
	\path (y3) ++(-120:\ty*\rr) coordinate [label=right:$y_7$] (y7);
	\path (y4) ++(180:\ty*\rrr) coordinate (y8);
	\path (y4) ++(60:\ty*\rrr) coordinate (y9);
	\path (y5) ++(60:\ty*\rrr) coordinate (y10);
	\path (y5) ++(-60:\ty*\rrr) coordinate (y11);
	\path (y6) ++(60:\ty*\rrr) coordinate (y12);
	\path (y6) ++(-60:\ty*\rrr) coordinate (y13);
	\path (y7) ++(-60:\ty*\rrr) coordinate (y14);
	\path (y7) ++(180:\ty*\rrr) coordinate (y15);

        \path (y8) ++(-120:\ty*\rrrr) coordinate (y16);
	\path (y8) ++(120:\ty*\rrrr) coordinate (y17);
	\path (y9) ++(120:\ty*\rrrr) coordinate (y18);
	\path (y9) ++(0:\ty*\rrrr) coordinate (y19);
	\path (y10) ++(120:\ty*\rrrr) coordinate (y20);
	\path (y10) ++(0:\ty*\rrrr) coordinate (y21);
	\path (y11) ++(0:\ty*\rrrr) coordinate (y22);
	\path (y11) ++(-120:\ty*\rrrr) coordinate (y23);
        \path (y12) ++(120:\ty*\rrrr) coordinate (y24);
	\path (y12) ++(0:\ty*\rrrr) coordinate (y25);
	\path (y13) ++(0:\ty*\rrrr) coordinate (y26);
	\path (y13) ++(-120:\ty*\rrrr) coordinate (y27);
	\path (y14) ++(0:\ty*\rrrr) coordinate (y28);
	\path (y14) ++(-120:\ty*\rrrr) coordinate (y29);
	\path (y15) ++(-120:\ty*\rrrr) coordinate (y30);
	\path (y15) ++(120:\ty*\rrrr) coordinate (y31);

  % fork:
  \def\fork{\path[draw, color=blue, ultra thick]}

    \fork (y2) -- (y1) -- (y3);
    \fork (y4) -- (y2) -- (y5);
    \fork (y6) -- (y3) -- (y7);

   \path[draw, color=blue, ultra thick] (y0) -- (y1);

  \foreach \n [evaluate=\n as \m using {int(2*\n)}] in {4,5,6,7} { 
    \path[draw, color=blue, ultra thick, dotted] (y\n) -- (y\m);
  }

 \foreach \n [evaluate=\n as \m using {int(2*\n+1)}] in {4,5,6,7} { 
    \path[draw, color=blue, ultra thick, dotted] (y\n) -- (y\m);
  }

  \fill [black] (y0) circle (2pt);

  \foreach \n in {1,2,...,15} {
    \fill [blue] (y\n) circle (2pt);
  }
  
  \foreach \n in {16,17,...,31} {
    \fill [black] (y\n) circle (2pt);
  }
% \end{scope}
 \end{tikzpicture}
  \end{center}
  \caption{Planar universal Steiner tree}
  \label{fig:bi3}
\end{figure}
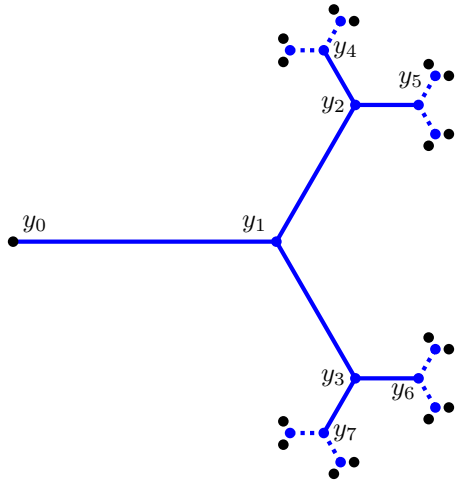

Thus almost any questions in high dimensions remain open. In particular, I would like to generalize the examples from our paper~\cite{cherkashin2025steiner}:

%I am interested in the generalization of results of~\cite{cherkashin2025steiner}:% in a solution to the Steiner problem for an input which is a mix of a regular $n$-gon with the setup of this paper. 
\begin{question}
For a regular $3$-dimensional $n$-gon in the plane $x = 1$ centered at $(1,0,0)$; let $\mathcal{Q}$ be the set of its vertices and consider the union of all $f^k(\mathcal{Q})$, $k \in \mathbb{N}$, where $f$ is a homothety with the center in the origin and a small enough scale factor $\lambda$. What is the structure of the corresponding Steiner tree?
   % Consider a regular $3$-dimensional $n$-gon in the plane $x = 1$ centered at $(1,0,0)$; let $\mathcal{Q}$ be the set of its vertices.
%Then the input is the union of all $f^k(\mathcal{Q})$, $k \in \mathbb{N}$, where $f$ is a homothety with the center in the origin and a small enough scale factor $\lambda$.
\end{question}

Also the generalization of universal Steiner tree (for planar ones see~\cite{paolini2015example},~\cite{cherkashin2023self}) is quite interesting. 
\begin{question}
   % To provide a nonplanar example of unique full Steiner tree with infinite number of branching points.
    Can one construct a non-planar example of a unique full Steiner tree with infinitely many branching points?
\end{question}
Work is underway on both issues mentioned above, and some progress has been achieved.
%Work is underway on both issues mentioned in this paragraph and there has been some progress achieved.
%\hline

%\noindent\makebox[\linewidth]{\rule{\paperwidth}{0.4pt}}
%\input{trash2}

%\paragraph{Acknowledgements.} 
%Y. Teplitskaya is supported from the French National Research Agency (ANR) under grant ANR-21-CE40-0013-01 (project GeMfaceT).
\paragraph{Acknowledgements.}

This research is supported by the French National Research Agency (ANR) under grant ANR-21-CE40-0013-01 (project GeMfaceT).
I would like to thank all colleagues mentioned in the paper for inspiring discussions, and especially Danila Cherkashin for kind assistance in preparing the manuscript.
\bibliography{main}
\bibliographystyle{plain}

\end{document}